\font\gothic=eufm10.
\font\piccolo=cmr8.
\font\script=eusm10.
\font\sets=msbm10.
\font\stampatello=cmcsc10.
\def\0{{\bf 0}}
\def\1{{\bf 1}}
\def\defineq{\buildrel{def}\over{=}}
\def\EqByDef{\buildrel{\bullet}\over{=}}

\def\QED{\hfill {\rm QED}}
\def\square{\hbox{\vrule\vbox{\hrule\phantom{s}\hrule}\vrule}}
\def\C{\hbox{\sets C}}
\def\N{\hbox{\sets N}}

\def\Primes{\hbox{\sets P}}

\def\Z{\hbox{\sets Z}}
\def\supporto{{\rm supp}}

\def\BH{(\hbox{\stampatello BH})}

\def\RamaExp{(\hbox{\stampatello Ramanujan Expansion})}
\def\LuchtExp{(\hbox{\stampatello Lucht Expansion})}

\def\DEE{(3.1)}
\def\Eq*{(4.0)}
\def\WN{(4.1)}
\def\UN{(4.2)}
\def\CI{(4.3)}
\def\TEQone{(5.1)}
\def\TEQtwo{(5.2)}
\def\TEQAone{(5.1)_A}
\def\TEQAtwo{(5.2)_A}
\def\TEQthree{(5.3)}
\def\TEQfour{(5.4)}
\def\HLzero{(2.1)}
\def\HLone{(6.1)}
\def\HLtwo{(6.2)}
\def\HLthree{(6.3)}
\def\HLfour{(6.4)}

\def\ODD{\hbox{\piccolo ODD}}
\def\Car{{\rm Car}}
\def\Win{{\rm Win}}
\def\CarT{\Car\; }
\def\WinT{\Win\; }

\def\IPP{(\hbox{\stampatello IPP})}

\par
\centerline{\bf A new class of Correlations}
\par
\centerline{\bf insisting on Ramanujan expansions}
\bigskip
\centerline{Giovanni Coppola}\footnote{ }{MSC $2020$: $11{\rm N}37$ - Keywords: correlation, Diophantine equation, periodicity, Ramanujan expansion} 

\bigskip

\par
\noindent
{\bf Abstract}. {\piccolo Studying Correlations with Ramanujan Expansions, we arrive to present the new class of, say, Two-Seasons Correlations, abbr.  T-S, as a natural set expressing some of the features of, say, H-L-like Correlations; these are the ones that mimic the H-L ($=$Hardy-Littlewood) Correlation with shift $2k$, needed to study $2k-$twin primes following Hardy \& Littlewood Conjecture. After introducing the $3-$Hypotheses Correlations in a previous paper, we add two other, very natural, hypotheses: the fifth is a technical one, simplifying calculations; but the fourth is called 'Parity', since it deals with the parity of natural numbers we play with. In particular, we may build (devoting to this 'our mainstream', here) a single Correlation that satisfies these '5 Axioms', thus a T-S one, that 'entangles two different Correlations' (whence Two-Seasons: T-S) depending on $a$ ($=$ the shift) parity. For $a$ even, our 'Artifact' mimics the H-L Correlation, in fact $a=2k$; but, while H-L Correlation is 'negligible', say, on $a$ odd, our Artifact seems to compare at least in the order of magnitude to H-L Correlation on $a$ even, being linked to another additive problem. Namely, on $a$ even, the Artifact 'counts', say, classic solutions to: $p_1+a=p_2$, in odd primes $p_1,p_2$; while, on $a$ odd, it 'counts' solutions to: $p_1+a=2^j p_2$, again with odd primes $p_1,p_2$ and with $j\in \N$ \hfill (satisfying the natural arithmetic constraints). \hfill More in general, our T-S Correlations 'entangle' two different Diophantine equations.}  

\bigskip
\bigskip
\bigskip

\par
\noindent{\bf 1. Introduction and general ideas. Notations, a new notation and its elementary properties} 
\bigskip
\par
\noindent
We may consider, once {\stampatello fixed} $f:\N \rightarrow \C$ and $g:\N \rightarrow \C$, namely two arbitrary {\stampatello arithmetic functions}, with a {\stampatello given parameter} $N\in \N$, their {\stampatello Correlation}, 
$$
C_{f,g}(N,a)\defineq \sum_{n\le N}f(n)g(n+a), 
\enspace \forall a\in\N, 
$$
\par
\noindent
see next section for details (and compare [C1]); our previous paper [C1] gave a short panorama of elementary methods, to study these, say \lq \lq pervasive\rq \rq, so to speak, objects of Analytic Number Theory. In applying a sprinkle of standard elementary number theory, esp., Dirichlet characters, we relied there upon (mainly) three hypotheses, for our Correlations: say, we were considering $3-$hp.s Correlations (assuming three hypotheses) in calculating, so-to-speak, the expected \lq \lq right Ramanujan Coefficients\rq \rq, good for example in obtaining at least heuristic formul\ae, for our $C_{f,g}(N,a)$ (the definition of Ramanujan Coefficients follows in a moment). As we saw in [C1], the first step is: assume that $g'$, the {\stampatello Eratosthenes Transform} [C1] of our $g$, vanishes after a given $Q\le N$; joining a technical requirement to our Correlation (namely, it has to be {\stampatello Fair}, see $\S2$), we start harvesting, in an elementary fashion, main properties of $C_{f,g_Q}(N,a)$, with $g_Q$ the $Q-${\stampatello Truncation} of our $g$, $\S2$, see Theorem 3.1 in [C1]. These two requirements on our Correlation (the $g-$Truncation \& Fairness) together are referred to as : the {\stampatello Basic Hypothesis}, abbrev. $\BH$, see [C1]. After Truncation, possible whatever is $g$, we get a great simplification : working with a Truncated Divisor Sum, compare $\S7$ results, is necessary to quoted Th.m 3.1 [C1], but also to present results in $\S4$, that link the periodicity of $g$ to the periodicity of $C_{f,g}(N,a)$, w.r.t. (with respect to) $a\in\N$. \hfill After $\BH$, [C1] adds hypotheses $(1),(2)$. 
\par
\noindent
This makes possible the calculation of {\it Wintner's Coefficients} (compare [C1],[C2],[C3]) for the IPP-{\it ification} of our $C_{f,g}(N,a)$, namely (compare [C1] definition) the function: \enspace $\widetilde{C_{f,g}}(N,a)\defineq C_{f,g}(N,\kappa(a))$,\enspace $\forall a\in\N,$ \enspace where $\kappa(n)\defineq \prod_{p|n}p=$ {\it square-free kernel} of any $n\in \N$, with $\kappa(1)\EqByDef 1$ as {\bf empty products are} $1$. We'll hereafter abbreviate $\EqByDef$ for \lq \lq equals by definition\rq \rq. Namely, we are quoting Corollary 3.1 [C1], holding under the three hypotheses $\BH$, $(1)$ \& $(2)$. The study of Correlations {\stampatello insists on Ramanujan Expansions} : RE.s, and not the Author, \lq \lq insist\rq \rq, say, in proposing themselves as the most powerful tool for Correlations. 
\par
We call $G:\N \rightarrow \C$ $\underline{\hbox{\stampatello a\enspace Ramanujan\enspace Coefficient}}$ for any given $F:\N \rightarrow \C$ whenever (sum $q\le x$, then send $x\to \infty$)
$$
\forall a\in \N, \quad
F(a)=\sum_{q=1}^{\infty}G(q)c_q(a)
$$
\par
\noindent
(given $F$, {\it uniqueness} of $G$ is {\it not guaranteed}: [C3]), with $c_q(a)\defineq \sum_{j\le q}^{\ast}\cos(2\pi ja/q)$ {\stampatello the Ramanujan Sum}. 

\vfill
\eject

\par				
\noindent
By the way, hereafter a sum with a $\ast$ above indicates {\it reduced residue classes}, w.r.t. the implicit module; here it's $q\in \N$ : in case of $c_q(a)$, the {\it argument} may be integer, $a\in \Z$. \hfill For Correlations, {\it shift} $a\in \N$, see $\S2$.
\par
We wrote \lq \lq {\stampatello a Ramanujan Coefficient}\rq \rq, since uniqueness holds very rarely: only assuming properties for $G$ (and not for $F$, that's not useful at all!), see [C0] and [C3]. 
\par
We specified, in the above definition of a Ramanujan Coefficient, the {\it classic method of summation}: the partial sums are over $q\le x$, then $x\to \infty$. This is NOT the only method of summation: compare the \lq \lq {\stampatello smooth summation method}\rq \rq, see [C2], summing over $q\in (P)$, where $(P)$ is the set of {\stampatello $P-$smooth numbers} for any fixed prime number $P\in \Primes$ (the primes' set, hereafter $p$ with/without subscripts is a prime)
$$
(P)\defineq \{ n\in \N : p|n\enspace \Rightarrow\enspace p\le P\}. 
\quad
\hbox{\rm Notice}:\enspace 1\in (P).
$$
\par
\noindent
See [C2], for the smooth summation (in which we send $P\to \infty$, within primes) and for any kind of summation: all giving the same result, in case of {\it absolute convergence}; but without, may get unexpected behaviors, of our series above: see [C2],$\S1$ for details. 
\hfill 
We tell about other ideas in next sections: see next Paper's Plan. 
\par
\centerline{\bf Notations, new notation with properties}
\par
\noindent
We use a standard Landau notation $O$, rather synonymous to Vinogradov notation $\ll$, both with subscripts indicating dependence on some parameters : see [T]; typically, depending on a positive arbitrarily small constant $\varepsilon$: using $O_{\varepsilon}$ or $\ll_{\varepsilon}$. With $\ast$ the standard Dirichlet convolution product and $\mu$ M\"obius function [T] we define the {\stampatello Eratosthenes Transform} (Wintner's terminology, we abbreviate E.T.) of a general $F:\N\rightarrow \C$ (an {\stampatello Arithmetic Function}) as $F'(d)\defineq \sum_{t|d}F(t)\mu(d/t)$, $\forall d\in \N$ ($d$ abbreviates {\stampatello divisor}). 
By the way, the {\stampatello divisors} have an important r\^ole for the \lq \lq shift-carrying factor\rq \rq, in Correlations, see $\S2$ and $\S3$.
\par
Another abbreviation we use is T.D.S. for {\stampatello Truncated Divisor Sum} : these are {\bf defined} in $\S7$ $\underline{\hbox{\stampatello Lemma}}$; $\0$ abbreviates {\stampatello the null-function}, i.e. $\0(n)\defineq 0$, $\forall n\in \N$. In $\S4$ we use l.c.m. for: least common multiple.   
\par
We define {\stampatello two operators}, strictly {\stampatello linked} together: 
$$
{\rm ODD} : n\in \N \mapsto n_{\rm ODD}\defineq n/2^{v_2(n)}\in \N
\qquad
\hbox{\rm among natural numbers}
$$
with $v_2(n)$ {\stampatello the $2-$adic valuation} of $n\in \N$
\par
\noindent
{\stampatello and}
$$
{}^{\ODD} : F\in \C^{\N} \mapsto F^{\ODD}\defineq (\1_{\rm ODD}\cdot F')\ast \1\in \C^{\N}
\qquad
\hbox{\rm among arithmetic functions}
$$
\par
\noindent
in fact, $\C^{\N}$ is the set of {\stampatello arithmetic functions}, hereafter $\1_{\rm ODD}$ is the {\stampatello indicator function of odd naturals}, with $\1_{\hbox{\script S}}$ indicating any {\script S}$\subseteq \N$ (and $\1=\1_{\N}$ the {\stampatello constant} $1$ {\stampatello function}, $F'$ the E.T. of our $F$, we saw above and $\ast$ the {\stampatello Dirichlet convolution product} [T]). Hereafter, $\cdot$ is the pointwise product in $\C^{\N}$. 
\par
\noindent
See that {\stampatello odd naturals are divisor-closed} ($\hbox{\script S}\subseteq \N$ is such, by definition, if $n\in \hbox{\script S},d|n$ imply $d\in \hbox{\script S}$), so:
$$
F^{\ODD}(n)\EqByDef \sum_{{d|n}\atop {d\,{\rm ODD}}}F'(d)
 =\sum_{d|n_{\rm ODD}}F'(d)\EqByDef F(n_{\rm ODD})
\quad
\forall n\in \N.
\leqno{(1.1)}
$$
\par
\noindent
Last but not least, {\stampatello recall} $)P(\defineq \{n\in \N : p|n \Rightarrow p>P\}$ the set of $P-${\stampatello sifted} naturals, whence ({\stampatello notice} $1\in )P($, here)\enspace we get \enspace  $\forall n\in \N,\quad n=n_{(P)}n_{)P(},\enspace \hbox{\stampatello with}\enspace n_{(P)}\in (P)\enspace \hbox{\stampatello and} \enspace n_{)P(}\in )P($ ; choosing $P=2$, $\forall n\in \N$,\enspace $n=n_{(2)}n_{)2(}=2^{v_2(n)}n_{\rm ODD}$ : so, actually, we may also write\enspace $n_{\rm ODD}=n_{)2(}$ , $\forall n\in \N$. 
\vfill
\par
\centerline{\bf Paper's Plan}
\smallskip
\item{$\ast$} next section describes classes of Correlations, starting from one Hypothesis only, arriving to the five hypothesis Class called T-S: Two-Seasons Correlations; 
\item{$\ast$} in $\S3$, we tell about two arithmetic factors, of a Correlation, which describes two Diophantine Equations; 
\item{$\ast$} in $\S4$ we offer another application of T-S, from their periodicity properties, to get combinatorial identities; 
\item{$\ast$} in $\S5$ we summarize our main results for T-S Correlations: these include the T-S called the {\stampatello Artifact};
\item{$\ast$} in $\S6$ we give Proofs for our main results for T-S; after, a panorama of Hardy-Littlewood models;
\item{$\ast$} in $\S7$ we provide main properties we need for the Truncated Divisor Sums. 
\eject

\par				
\noindent{\bf 2. Classes of Correlations: Basic Hypothesis, $3-$Hypotheses, Two-Seasons} 
\bigskip
\par
\noindent
We recall the definition of {\stampatello Correlation}, its {\stampatello shift} $a\in \N$ (also called {\it argument}), its {\stampatello length} $N\in \N$, its {\stampatello two factors}, namely $f:\N \rightarrow \C$ {\stampatello and} $g:\N \rightarrow \C$, whose Correlation is in fact
$$
C_{f,g}(N,a)\EqByDef \sum_{n\le N}f(n)g(n+a), 
\enspace \forall a\in\N, 
$$
\par
\noindent
also called {\stampatello Shifted Convolution Sum} (compare [C1]), of our $f$ \& $g$; here we follow the [C1] definitions, also given \& recalled in $\S4$, following (see $\S4.1$, esp.). However, {\stampatello we recall BH}, the {\stampatello Basic Hypothesis}: we abbreviate with \enspace \lq \lq $g$\enspace {\stampatello has range}\enspace $Q\le N$\rq \rq, hereafter, \lq \lq $\supporto(g')\subseteq [1,Q] ,\thinspace Q\le N$\rq \rq, 
$$
g\enspace \hbox{\stampatello has\enspace range}\enspace Q\le N
\quad
\hbox{\stampatello and}
\quad
C_{f,g_Q}(N,a)=\sum_{q\le Q}\widehat{g_Q}(q)\sum_{n\le N}f(n)c_q(n+a)\enspace \hbox{\stampatello is\enspace Fair}
\leqno{\BH}
$$
\par
\noindent
which means, {\stampatello by definition: in equation above} (i.e., $(\ast)$ in [C1], $\S4.1$) {\stampatello the $a-$dependence is} $\underline{\hbox{\stampatello only}}$ {\stampatello in } $c_q$\enspace {\stampatello argument, i.e.: inside} $c_q(n+a)$ (in other words, $a-$dependence in the above is allowed nor in terms of $f(n)$, $\widehat{g_Q}(q)$, neither of their supports). Details about $\BH$ may be found in [C1]. 
\par
Whenever $g$ is not of range $Q\le N$, {\stampatello Theorem 2.1 [C1] (The Divisors' Cut)} allows to pass {\stampatello from } $g:\N \rightarrow \C$, a general arithmetic function, {\stampatello to its $N-$truncated divisor sum, called } $g_N$ , [C1]; and this {\stampatello passage is \lq \lq at small cost\rq \rq:} after cutting divisors, the error's small.
\smallskip
\par
{\stampatello Thus,} the {\stampatello real problem} is {\stampatello to have a Fair Correlation}.
\smallskip
\par
\noindent
In [C1], {\stampatello All \lq \lq Basic Hypothesis\rq \rq-Correlations'} main properties are {\stampatello proved in Theorem 3.1}. It {\stampatello supplies}, in particular, an {\stampatello explicit formula for} the {\stampatello Carmichael-Wintner coefficients} of our Correlation: see point $C)$ of quoted Th.m. 
\par
In the following, we'll indicate $\BH$ as {\stampatello Hypothesis } $(0)$, writing $g_Q$ for $g$ {\stampatello of range } $Q$ (here $Q\le N$); then, calling an arithmetic function $F$ $\IPP$, by def., when $F(a)=F(\kappa(a))$, $\forall a\in \N$ : see [C1,[C2], we have 
\par
\noindent
{\stampatello Hypothesis} $(1)$: $g_Q$ $\IPP$ , {\stampatello equivalent} (see $\S4.2$ beginning in [C1]) {\stampatello to both:}
$$
g'_Q=\mu^2\, \cdot\, g'_Q
\quad \Longleftrightarrow \quad
\widehat{g_Q}=\mu^2\, \cdot\, \widehat{g_Q}
$$
\par
\noindent
{\stampatello and}
\par
\noindent
{\stampatello Hypothesis} $(2)$: $f=\1_{\Primes}\,\cdot\, f$\quad  (namely, $f$ is {\stampatello supported on primes})

\medskip

\par
These \enspace $(0),(1),(2)$ \enspace {\stampatello are the requirements defining the} \lq \lq {\stampatello Three Hypotheses Correlations}\rq \rq. 

\medskip

\par
\noindent
The $3-hp.s-${\stampatello Correlations}, as we'll {\stampatello abbreviate those with $(0),(1)$ \& $(2)$,} are the object of study of [C1], for {\stampatello two main reasons : first, they} somehow {\stampatello mimic the Hardy-Littlewood Correlation}:
$$
C_{\Lambda,\Lambda}(N,a)\defineq \sum_{n\le N}\Lambda(n)\Lambda(n+a),
\enspace \forall a\in \N, \enspace
\hbox{\stampatello say,\enspace H-L Correlation}, 
\leqno{\HLzero}
$$
\par
\noindent
where von Mangoldt function is \thinspace $\Lambda(n)\defineq \1_{\Primes}\log \kappa(n)$, $\forall n\in \N$\thinspace [T]; we wish to study this H-L, when following $a=2K$ for $2K-$twin primes count, given {\stampatello by Hardy \& Littlewood in } 1923 [HL], {\stampatello Conjecture B}.
\par
{\stampatello Second} reason, {\stampatello a technical issue}: 
$$
C_{\Lambda,\Lambda}(N,a)\sim \hbox{\gothic S}(a)N , 
\enspace \hbox{\stampatello as } N\to \infty,\enspace \forall a\in \N
\enspace \hbox{\stampatello fixed\thinspace and\thinspace even,}\enspace \hbox{\rm is\thinspace equivalent\thinspace to}\enspace \hbox{\stampatello Conjecture\enspace B,\thinspace with}
$$
$$
\hbox{\gothic S}(a)\defineq \sum_{q=1}^{\infty}{{\mu^2(q)}\over {\varphi^2(q)}}c_q(a),\enspace \forall a\in \N
$$
\par
\noindent
the {\stampatello classic Singular Series}; since $q$ here is square-free, we have (see Fact 3.2 in [C2] for an elementary Proof) $\forall a\in \N$ fixed, $c_q(a)=c_q(\kappa(a))$, whence the singular series $\hbox{\gothic S}(a)$ $\IPP$, i.e., $\hbox{\gothic S}(a)=\hbox{\gothic S}(\kappa(a))$: it \lq \lq Ignores Prime Powers\rq \rq, inside $a$, $\forall a\in \N$.
\par				
\noindent
In [C1] the $3-hp.s-${\stampatello Correlations} $C_{f,g_Q}(N,a)$ {\stampatello have \lq \lq IPP-ification\rq \rq\enspace} $C_{f,g_Q}(N,\kappa(a))$ {\stampatello that has same Wintner's coefficients of original} $C_{f,g_Q}(N,a)$, as we prove in {\stampatello Corollary 3.1 [C1].} In some sense, it is {\stampatello a kind of expected behaviour, since Wintner's coefficients of} $C_{\Lambda,\Lambda}(N,a)$ {\stampatello \lq \lq build up\rq \rq\enspace main term} $\hbox{\gothic S}(a)N$ in the {\stampatello H-L Conjecture}.

\bigskip

\par
Now, {\stampatello passing from $3-hp.s-$Correlations to T-S Correlations, we add to $(0)$, $(1)$, $(2)$ only two further hypotheses.} However, with \lq \lq {\stampatello highly unexpected}\rq \rq \thinspace {\stampatello effects}, {\stampatello as we'll see} in a moment. 

\medskip

\par
\noindent
Since hypothesis $(3)$ in [C1] {\stampatello regarded Ramanujan Conjecture}, while our present fourth hypothesis regards {\stampatello parity \& not the above}, we'll {\stampatello change terminology}. We explicitly, here, mean by Axiom a hypothesis. 

\medskip

\par
We'll {\stampatello speak about Axioms. We list them as $(0)$, $(1)$ \& $(2)$ like above hypotheses.} However, $(3)$ \& $(4)$ {\stampatello are our fourth \& fifth, adding properties to our first three.} 
\par
Last but not least : $(0)$, $(1)$ \& $(3)$ {\stampatello Axioms contain restrictions on} $g'$ {\stampatello (not on $\widehat{g}$, even if we'll see the links)}, $(2)$ \& $(3)$ {\stampatello Axioms give constraints to } $f$ ; $(4)$ {\stampatello completes the list, with a simplifying technical requirement, also involving } $(0)$ (so to \lq \lq {\stampatello close the circle}\rq \rq). 

\bigskip

\par
\noindent
{\stampatello Axiom} $(0)$. \enspace $\supporto(g')\subseteq [1,Q]$ , $Q\le N$ \enspace {\stampatello and } $C_{f,g}(N,a)$ {\stampatello is Fair} 
\par
\noindent
{\stampatello Axiom} $(1)$. \enspace $g'=\mu^2\,\cdot\,g'$ (i.e., $g'$ {\stampatello is supported within square-free naturals})
\par
\noindent
{\stampatello Axiom} $(2)$. \enspace $f=\1_{\Primes}\,\cdot\,f$ (i.e., $f$ {\stampatello is supported within prime numbers})
\par
\noindent
{\stampatello Axiom} $(3)$. \enspace $g'=\1_{\ODD}\,\cdot\,g'$ (i.e., $g'$ {\stampatello is odd-supported}) \enspace {\stampatello and} \enspace  $f=\1_{\ODD}\,\cdot\,f$ (i.e., $f$ {\stampatello is odd-supported})
\par
\noindent
{\stampatello Axiom} $(4)$. {\stampatello Above $Q\in \N$ in $(0)$ is } $Q=N\ge 9$;\enspace {\stampatello and both } $N\not\in \Primes$, $N-1\not\in \Primes$

\medskip

\par
{\stampatello Obviously we call Axiom $(3)$ (Parity)}. Axioms $(0)$ {\stampatello to} $(3)$ {\stampatello regard} $g'$ \& $f$ {\stampatello supports}, including the requirement \lq \lq {\stampatello Fair}\rq \rq, implicitly depending on $\widehat{g}$ \& $f$ (so, on $g'$ \& $f$); while {\stampatello fourth's for parameters}. 

\medskip

\par
\noindent
{\stampatello Summarizing, all} $\underline{\hbox{\stampatello Correlations\enspace with\enspace Axioms\enspace} (0),(1),(2),(3),(4)}$ {\stampatello will be} $\underline{\hbox{\stampatello called\enspace Two-Seasons}}$ {\stampatello Correlations}, which will be {\stampatello abbreviated as T-S}.
\smallskip
\par
Of course, {\stampatello both Theorem 3.1 \& Corollary 3.1 in } [C1] {\stampatello still hold for T-S}. 

\medskip

\par
\noindent
We may reformulate these Axioms saying that : $C_{f,g}(N,a)$ {\stampatello is\enspace T-S}, by definition, if and only if
$$
C_{f,g}(N,a) \thinspace \hbox{\stampatello Fair}; \thinspace \supporto(g')\subseteq \{d\le N : \mu^2(d)=1,d\;\, \ODD\}; \thinspace \supporto(f)\subseteq \{p\in \Primes : p\;\, \ODD\}; N\ge 9; N,N-1\not\in \Primes 
$$

\medskip

\par
Also, we may give the same restrictions, in {\stampatello Axioms} $(0)$, $(1)$ \& $(3)$, assumed on $g'$ but for $\widehat{g}$; i.e., 
\medskip
\par
\noindent
\centerline{in $(0)$: \enspace $\widehat{g}=\1_{[1,Q]}\cdot \widehat{g}$,\hfill in $(1)$: \enspace $\widehat{g}=\mu^2\cdot \widehat{g}$\hfill and in $(3)$: $\widehat{g}=\1_{\ODD}\cdot \widehat{g}$.} 
\smallskip
\par
\noindent
Hereafter we abuse notation for indicator functions, as we'll write $\1_{\hbox{\script R}}$ to abbreviate $\1_{\hbox{\script R}\cap \N}$, whatever is the subset of real numbers {\script R} : usually, an interval (here closed). In other words $\1_{[1,Q]}$ indicates the set $\{1,\ldots,Q\}$. 
\medskip
\par
This is a kind of \lq \lq {\stampatello switch}\rq \rq \thinspace in between $g'$ and $\widehat{g}$ which {\stampatello follows from}: $g$ is a {\stampatello T.D.S.} ({\stampatello Truncated Divisor Sum}); in fact, under this assumption on $g$, our $\underline{\hbox{\stampatello Lemma}}$ in the $\S7$, {\stampatello Appendix}, proves the {\stampatello switch}: for details, see the considerations in $\S4.1$, soon after {\stampatello Wintner's Period} definition. 
\medskip
\par
(Furthermore, even if we haven't proved the implication: $\widehat{g}=\1_{[1,Q]}\cdot \widehat{g}$ $\Rightarrow$ $g'=\1_{[1,Q]}\cdot g'$, in $\S7$ $\underline{\hbox{\stampatello Lemma}}$, same arguments there, following the $\LuchtExp$ in $\S7$ $\underline{\hbox{\stampatello Lemma}}$ Proof, show it easily.) 

\vfill
\eject

\par				
\noindent{\bf 3. Two-Seasons Correlations' Axioms entangle two Diophantine equations} 
\bigskip
\par
\noindent
We start with two {\stampatello easy arithmetic functions}:
$$
f=\1_{\hbox{\script F}}\cdot \1_{\rm ODD}
\qquad \hbox{\stampatello and} \qquad
g=\1_{\hbox{\script G}}^{\ODD} ,
\leqno{\DEE}
$$
\par
\noindent
where \enspace $\1_{\hbox{\script F}}$ \enspace is the {\stampatello indicator of a fixed}\enspace {\script F}$\subseteq \N$; like, instead, \enspace $\1_{\hbox{\script G}}$ \enspace {\stampatello indicates a fixed}\enspace {\script G}$\subseteq \N$.
\par
However, we {\stampatello consider} \enspace $\1_{\hbox{\script G}}^{\ODD}(m)\EqByDef \1_{\hbox{\script G}}(m_{\rm ODD})\EqByDef \1_{\hbox{\script G}}(m/2^{v_2(m)})$, $\forall m\in \N$. 
\par
\noindent
{\stampatello Thus}
$$
C_{f,g}(N,a)=\sum_{{n\le N}\atop {{n\,{\rm ODD}}\atop {n\in \hbox{\script F}}}}\sum_{0\le j\le {{\log(N+a)}\over {\log 2}}}\1_{\rm ODD}\left({{n+a}\over {2^j}}\right)\1_{\hbox{\script G}}\left({{n+a}\over {2^j}}\right)
 =\sum_{0\le j\le {{\log(N+a)}\over {\log 2}}}\sum_{{n\le N, n\,{\rm ODD}, n\in \hbox{\script F}}\atop {{{n+a}\over {2^j}}\,{\rm ODD}, {{n+a}\over {2^j}}\in \hbox{\script G}}}1,
\enspace \forall a\in \N
$$
\par
\noindent
whence ($a$ {\stampatello even} $\Rightarrow$ $j=0$) : 
$$
C_{f,g}(N,a)=\sum_{{n\le N, n\,{\rm ODD}, n\in \hbox{\script F}}\atop {n+a\in \hbox{\script G}}}1,
\qquad \quad \forall a\in \N\enspace \hbox{\stampatello even}
\eqno{(\hbox{\stampatello regarding\thinspace \thinspace a\thinspace \thinspace \lq \lq single\rq \rq \thinspace \thinspace Diophantine \thinspace \thinspace Equation})}
$$
\par
\noindent
and ($a$ {\stampatello odd} $\Rightarrow$ $1\le j\le \log(N+a)/\log 2$)
$$
C_{f,g}(N,a)=\sum_{1\le j\le {{\log(N+a)}\over {\log 2}}}\sum_{{n\le N, n\,{\rm ODD}, n\in \hbox{\script F}}\atop {{{n+a}\over {2^j}}\,{\rm ODD}, {{n+a}\over {2^j}}\in \hbox{\script G}}}1,
\quad \forall a\in \N\enspace \hbox{\stampatello odd}
\eqno{(\hbox{\stampatello regarding\thinspace \thinspace \lq \lq one\rq \rq \thinspace \thinspace different\thinspace \thinspace D. E.})}
$$

\medskip

\par
\noindent
The functions in $\DEE$, {\stampatello actually}, {\stampatello satisfy} {\stampatello Parity Axiom} (the \lq \lq third\rq \rq!). 

\medskip

\par
{\stampatello Thus}
\par
\noindent
{\stampatello Diophantine Entanglement}, so-to-speak, {\stampatello starts from Parity}.

\bigskip

\par
\noindent
We gave a simple example of Correlation, here, just to give a crystal clear idea of how Axiom $3$ works; with other Axioms interaction, say, things get more complicated especially when coming at truncations, for the divisors (namely, Axiom $0$).  

\bigskip

\par
See $\S5.1$ for T-S Correlations having a similar behavior.  

\vfill
\eject

\par				
\noindent{\bf 4. Periodicity properties for Correlations' shifts: Wintner's Period. Combinatorial Identities} 
\bigskip
\par
\noindent{\bf 4.1. Correlations and Wintner's Period}
\smallskip
\par
\noindent
We recall (see $\S2$) that the Correlation of two fixed\enspace $f,g:\N \rightarrow \C$\enspace , with fixed length \enspace $N\in \N$ (we think about it as a large parameter, eventually $N\to \infty$), is \enspace $C_{f,g}(N,a)\EqByDef {\displaystyle \sum_{n\le N} }f(n)g(n+a)$,\enspace $\forall a\in \N,$ \enspace itself an arithmetic function, depending on the shift \enspace $a\in \N$; obviously, we {\stampatello call $g$ the shift-carrying factor} of \enspace $C_{f,g}(N,\bullet):a\in \N \mapsto C_{f,g}(N,a)\in \C$, while the other factor is $f$. 
\par
We soon study the {\stampatello periodicity} of \enspace $C_{f,g}(N,a)$, w.r.t. $a\in \N$, which is obviously strictly {\stampatello linked to} the {\stampatello periodicity of} \enspace $g(m)$, w.r.t. $m\in \N$. And the name {\stampatello shift-carrying factor} for  $g$ is not by chance.
\par
A beautiful idea comes into play, right now : since a {\stampatello Ramanujan Expansion} for \enspace $C_{f,g}(N,a)$, w.r.t. $a\in \N$, is not so easy to prove, at least for \lq \lq manageable\rq \rq, say, coefficients (compare Theorem 3.1 in [C1], that calculates the \lq \lq best candidates\rq \rq, namely the {\stampatello Carmichael-Wintner} coeff.s, \lq \lq for {\stampatello BH-}Correlations\rq \rq), we may apply the idea implicitly used in [CM], compare {\stampatello The Divisors' Cut} [C1] (Th.m 2.1 there); i.e., we pass {\stampatello from} a very {\stampatello general} \enspace $g:\N \rightarrow \C$ \enspace {\stampatello to its Truncated counterpart} \enspace $g_N:\N \rightarrow \C$, \enspace having no divisors beyond $N$ ($=$Correlation's Length), after their \lq \lq cut\rq \rq. 
\par
In the {\stampatello Appendix} (see $\S7$) we give, in fact, for a generic T.D.S., like $g_N$ here, a {\stampatello Ramanujan Expansion with fixed length}; these are not only {\stampatello finite Ramanujan Expansions} [CM], but actually they have {\stampatello Carmichael-Wintner coefficients} (supported within \enspace $[1,N]\cap \N$, for $g_N$ case): compare the $\underline{\hbox{\stampatello Lemma}}$ in $\S7$ and its {\stampatello Proof}, using two fixed-length Expansions ({\stampatello Ramanujan's \& Lucht's}).
\smallskip
\par
Here we need, for \enspace $g_N(m)\defineq {\displaystyle \sum_{{d\le N}\atop {d|m}} }g'(d)$, where \thinspace $g'$ \thinspace is {\stampatello $g$ Eratosthenes Transform}, with {\stampatello Ramanujan Coefficient} \enspace $\widehat{g_N}=\CarT\, g_N=\WinT\, g_N$, see quoted $\underline{\hbox{\stampatello Lemma}}$, the 
$$
\forall m\in \N, \quad
g_N(m)=\sum_{q\le N}\widehat{g_N}(q)c_q(m).
\leqno{\RamaExp}
$$
\par
\noindent
{\stampatello For any $g_N\neq \0$ (see $\S7$) we define} \enspace $W_Q(g_N)\defineq $l.c.m.$\{ q\le Q : \widehat{g_N}(q)\neq 0\}$ {\stampatello to be} $\underline{\hbox{\stampatello Wintner's Period}}$. 
\par
\noindent
In particular, $g_N\neq \0$ {\stampatello is constant} IFF(If and only If) $\widehat{g_N}(q)=0, \forall q>1$, IFF $W_N(g_N)=1$ (as $\widehat{g_N}(1)\neq 0$). 
\par
\noindent
Our definition's parameter \hfill $Q\in \N$ \hfill may be \hfill $Q\neq N$; anyway, \hfill $W_Q(g)$ \hfill is well-defined IFF \hfill $g\neq \0$ \hfill is a T.D.S. 
\par
The name comes from previous $g_N$ expansion, observing that: $\forall m\in \N$, $g_N(m+W_N(g_N))=g_N(m)$, by the {\stampatello $q-$periodicity of} $c_q(m)$, w.r.t. $m\in \N$. \hfill As we discover soon, this {\stampatello property's precious}. 
\par
For example, expanding \enspace $g_N$ \enspace as above entails (the same $(\ast)$ of [C1], $\S4.1$, i.e.) 
$$
C_{f,g_N}(N,a)=\sum_{q\le N}\widehat{g_N}(q)\sum_{n\le N}f(n)c_q(n+a),
\quad
\forall a\in \N. 
\leqno{\Eq*}
$$
\par
\noindent
{\stampatello Thus}
$$
C_{f,g_N}(N,a)=C_{f,g_N}(N,a+W_N(g_N)),
\quad
\forall a\in \N. 
$$
\par
\noindent
The possibility to have, here, a period \thinspace $<W_N(g_N)$ \thinspace is not excluded. 
\par
\noindent
In fact, {\stampatello it may happen} that {\stampatello we have a} \thinspace $f$ \thinspace {\stampatello depending on $W_N(g_N)$} and $a\in \N$, say \enspace $f=f_{a,N}$ : by $\Eq*$, 
$$
C_{f_{a,N},g_N}(N,a)=\sum_{q\le N}\widehat{g_N}(q)\sum_{n\le N}f_{a,N}(n)c_q(n+a)
 \defineq \sum_{q\le N}\widehat{g_N}(q)\sum_{{n\le N}\atop {n+a\equiv 0(\!\!\bmod W_N(g_N))}}f(n)c_q(n+a)
  :=\Sigma,
$$ 
\par
\noindent
say, choosing $f_{a,N}$ in this way {\stampatello by definition} ; whence, in turn, $W_N(g_N)$ definition entails that above's 
$$
\Sigma=\sum_{q\le N,\widehat{g_N}(q)\neq 0}\widehat{g_N}(q)\sum_{{n\le N}\atop {n+a\equiv 0(\!\!\bmod W_N(g_N))}}f(n)c_q(n+a)
 =\sum_{q\le N,\widehat{g_N}(q)\neq 0}\widehat{g_N}(q)\sum_{n\le N}f(n)\varphi(q), 
$$
\par
\noindent
that's {\stampatello constant} w.r.t. {\stampatello the shift $a\in \N$}; here, second equation {\stampatello uses in fact}:
$$
\left.W_N(g_N))\enspace \right|\enspace m
\quad \Rightarrow \quad
c_q(m)=\varphi(q), 
\enspace \forall q\in \supporto(\widehat{g_N}).
$$
\par				
\noindent
Hence, our Correlation has period $1$:\enspace $C_{f_{a,N},g_N}(N,a)=C_{f_{a,N},g_N}(N,a+1)$, $\forall a\in \N$, {\stampatello for $f=f_{a,N}$.} If we choose $g_N\neq \0$ not a constant, we have\enspace $W_N(g_N)>1$, whence our\enspace $C_{f,g_N}(N,\bullet)$\enspace has minimal period $1<W_N(g_N)$, at least for $f$ here. 
\smallskip
\par
(This example for $f$ is just, so to speak, {\stampatello \lq \lq pathologic\rq \rq, inasmuch $f$ depends on $a\in \N$}.) 

\medskip

\par
\noindent
{\stampatello At this point} we have to {\stampatello remark that \enspace $W_N(g_N)$ \enspace isn't, in general, the minimal period for } $g_N$. However, $g_N$ is {\stampatello $W_N(g_N)-$periodic}, compare soon after \enspace $W_N$ \enspace definition.
\par
\noindent
Also, {\stampatello our \enspace $W_N(g_N)$ \enspace is not, in general, the minimal period for } $C_{f,g_N}(N,\bullet)$ (namely, for the {\stampatello arithmetic function} sending $a\in \N$ to $a-${\stampatello Correlation}). 
\par
In all, {\stampatello without any statement about minimal periods}, we have: 
$$
\forall m\in \N, \enspace g_N(m)=g_N(m+W_N(g_N))
\quad \hbox{\stampatello and} \quad
\forall a\in \N, \enspace C_{f,g_N}(N,a)=C_{f,g_N}(N,a+W_N(g_N)). 
\leqno{\WN}
$$
\par
\noindent
Up to now, the {\stampatello only Axiom we used} is the Basic Hypothesis: $\BH$; actually, we used only the assumption that \enspace $g=g_N$, \enspace in fact our \enspace $f_{a,N}$ \enspace is {\stampatello excluded by our assumption} that the {\stampatello Correlation is Fair} : for $\BH$ and all the deeper aspects, see [C1], esp. its $\S2$. First equation in $\WN$, of course, is only the definition of $W_N(g_N)$ (using $g_N$ fixed length {\stampatello Ramanujan Expansion}), while the second in $\WN$, actually, needs \lq \lq something more\rq \rq, as it starts to be deeper: we'll use it, here and in future publications, under $\BH$.  

\medskip

\par
We look at {\stampatello Wintner's Period} and, hereafter, we {\stampatello focus on\enspace $W_N(g_N)$\enspace :} we leave all the deepening about \enspace $W_N(g_N)$\enspace \& its {\stampatello minimality} for future publications. 

\medskip

\par
It's possible to {\stampatello define } $W_N(g_N)$\enspace {\stampatello as soon as $g_N\neq \0$\enspace is a T.D.S.} We got \thinspace $g_N$ \thinspace from a general \thinspace $g$, \thinspace {\stampatello truncating} the {\stampatello E.T. (Eratosthenes Transform): } $g'$\enspace {\stampatello becomes } $g'\cdot \1_{[1,N]}$ (here, since\enspace $g'$\enspace is {\stampatello defined in} $\N$, we don't need to write $\1_{[1,N]\cap \N}$, instead).

\medskip

\par
\noindent
{\stampatello Axiom 0 , better a part of it,} namely a part of BH, {\stampatello multiplies } $g'$\enspace {\stampatello by } $\1_{[1,N]}$. 
\par
{\stampatello We add Axioms one by one.} 
\par
\noindent
{\stampatello Axiom 1 , i.e., $g \enspace \IPP$, multiplies, again } $g'$,\enspace {\stampatello by } $\mu^2$. 
\par
{\stampatello Thus first two Axioms have, for } $g'$,\enspace {\stampatello a support contained in square-free divisors} $\le N$. 
\par
\noindent
{\stampatello Axiom 3 , in particular, multiplies\enspace $g'$\enspace again, by } $\1_{\ODD}$. 
\par
{\stampatello These three Axioms make up a\enspace $g$\enspace with E.T.\enspace $g'$\enspace supported on odd \& square-free} $d\le N$. 

\medskip

\par
\noindent
The three sets : \enspace $\{ d\in \N : d\le N\}$, \enspace $\{ d\in \N : \mu^2(d)=1\}$, \enspace $\{ d\in \N : d\enspace \ODD \}$ \enspace {\stampatello are divisor-closed}, so : 
$$
g' \enspace \hbox{\stampatello has\enspace support\enspace in\enspace their\enspace intersection}
\quad \Rightarrow \quad 
\widehat{g} \enspace \hbox{\stampatello has\enspace support\enspace in\enspace their\enspace intersection}
$$
\par
\noindent
{\stampatello by} $\S7$ $\underline{\hbox{\stampatello Lemma}}$ \enspace $(\ast)$ \enspace \& \enspace $(\ast\ast)$; {\stampatello same result's } $(\ast)$ \enspace \& \enspace $(\ast\ast\ast)$ \enspace {\stampatello give converse,} \enspace \lq \lq $\Leftarrow$\rq \rq,\enspace {\stampatello too.}
\par
In this way , {\stampatello Axioms $(0)$, $(1)$ \& $(3)$} , {\stampatello in particular, ensure that}
$$
g_N\neq \0
\quad \Rightarrow \quad 
p\le N, \forall \left. p\;\right|\,W_N(g_N)\enspace \hbox{\stampatello and}\enspace W_N(g_N)\enspace \hbox{\stampatello is\enspace odd\enspace \& \enspace square-free} .
$$
\par
\noindent
In other words, the three {\stampatello Axioms} entail: 
$$
g_N\neq \0
\quad \Rightarrow \quad 
W_N(g_N)=\prod_{2<p\le N,p|q\in \supporto\left(\widehat{g_N}\right)}p, 
$$
\par
\noindent
abbreviating\enspace \lq \lq $p|q \in \supporto\left(\widehat{g_N}\right)$\rq \rq \enspace for the {\stampatello Condition}: \lq \lq $\exists q\in \supporto\left(\widehat{g_N}\right)$ \thinspace such \thinspace that $q\equiv 0(\bmod\enspace p)$\rq \rq. 
\par
More in general, the {\stampatello formula above entails} an easier periodicity, with $N-${\stampatello Universal period}: 
$$
g_N\neq \0
\quad \Rightarrow \quad 
\left.W(g_N)\,\right|\, U_N,
\quad
\hbox{\rm where} \enspace U_N\defineq \prod_{2<p\le N}p ;  
$$
\par
\noindent
whence,in particular, we have : 
$$
g_N\neq \0 ,\enspace \supporto(g_N)\subseteq \{d\in \N : d\;\ODD, \mu^2(d)=1\}
\quad \Rightarrow \quad 
\forall a\in \N, \enspace C_{f,g_N}(N,a)=C_{f,g_N}(N,a+U_N). 
\leqno{\UN}
$$
\par
\noindent
{\stampatello This} comes from: $\WN$ {\stampatello formul\ae \enspace and above property for } $W_N(g_N)$ \& $U_N$. And it {\stampatello sounds strange, as }$a$ \& $a+U_N$ {\stampatello have different parity}, $\forall a\in \N$. Also, $\UN$ $\Rightarrow$ $\underline{\hbox{\stampatello T-S\enspace Correlations\enspace are}}\enspace U_N-\underline{\hbox{\stampatello periodic}}$. 

\vfill
\eject

\par				
\noindent{\bf 4.2. Combinatorial Identities}
\bigskip
\par
\noindent
We'll come soon to see the Identities, which we call \lq \lq Combinatorial\rq \rq, since they look like their \lq \lq prototype\rq \rq, i.e. Equation $\UN$ in between Correlation's values with {\stampatello shifts} that {\stampatello are $U_N-$apart}.

\medskip

\par
We wish to give, now, some more information about the primes which certainly divide\enspace $W_N(g_N)$. This is a consequence of a {\stampatello very general property} of {\stampatello finite Ramanujan Expansions}, which we found in [CM] (see $\S3$, $(12)$ in particular) : it {\stampatello works for $g_N$ as} 
$$
\widehat{g_N}(q)={{g'(q)}\over q} , 
\quad 
\forall q\in \left({N\over 2}, N\right] ; 
$$
\par
\noindent
in fact, {\stampatello compare $\S2$ in [C1], these coefficients are Wintner's} 
$$
\widehat{g_N}(q)\defineq \sum_{{d\le N}\atop {d\equiv 0(\!\!\bmod q)}}{{g'(d)}\over d}
 = {1\over q}\sum_{m\le {N\over q}}{{g'(qm)}\over m}
$$
\par
\noindent
which, $\forall q>N/2$, have $m=1$, proving above formula for \enspace $N/2<q\le N$. Notice: these $q$ are natural numbers. 
\smallskip
\par
{\stampatello Thus, in the range}\enspace ${N\over 2}<q\le N$, \enspace $\widehat{g_N}(q)\neq 0$ \enspace {\stampatello is equivalent to:} \enspace $g'(q)\neq 0$.   
\smallskip
\par
\noindent
In particular, {\stampatello if\enspace $g_N\neq \0$, the Axioms $(0)$, $(1)$ \& $(3)$ entail} 
$$
\prod_{N/2<p\le N,\; g'(p)\neq 0}p=\left.\prod_{N/2<p\le N,\; \widehat{g_N}(p)\neq 0}p\enspace \right|\enspace W_N(g_N) , 
\qquad
\hbox{\stampatello as} \enspace N\ge 9. 
$$

\medskip

\par
\noindent
For example, our {\stampatello Artifact} has \enspace $g'(q):=-\mu(q)\log q$, so \enspace $g'(p)=\log p\neq 0$, \enspace $\forall p\in \left({N\over 2}, N\right]$. 

\medskip

\par
See that, say, 
$$
\left.\prod_{{N\over 2}<p\le N}p\enspace \right|\enspace W_N(g_N)
\quad \Rightarrow \quad
W_N(g_N) \gg \prod_{{N\over 2}<p\le N}p=e^{\theta(N)-\theta(N/2)} \gg e^{N/2}, 
$$
\par
\noindent
by the {\stampatello Prime Number Theorem (PNT)} in the form :\enspace $\theta(N)\sim N$, as \thinspace $N\to \infty$; see {\bf Fact 1.6} [C1].   

\bigskip

\par
\noindent
We wish, now, to get {\stampatello more interesting Combinatorial Identities,} esp. {\stampatello for T-S Correlations.} 

\medskip

\par
We'll {\stampatello confine to shifts} \enspace $a=1$ \& $a=2$, {\stampatello but using } $\underline{\hbox{\stampatello all}}$ {\stampatello our Axioms for T-S Correlations} will {\stampatello give links between \lq \lq not truncated\rq \rq} $g$, {\stampatello having \lq \lq small\rq \rq\enspace shifts } $a=1,2$\enspace {\stampatello and \lq \lq truncated\rq \rq} $g_N$, {\stampatello regarding \lq \lq huge\rq \rq\thinspace shifts } $a=U_N+1,U_N+2$. 

\medskip

\par
Our {\stampatello preliminary assumption avoids trivial cases } $g=\0$\enspace {\stampatello and } $g_N=\0$; {\stampatello then, we add the T-S Axioms} : {\stampatello in particular, } $(0)$, $(1)$ \& $(3)$ {\stampatello are used as hypotheses in } $\UN$, {\stampatello while } $(2)$ \& $(4)$ {\stampatello together, once applied in } $\S7$ $\underline{\hbox{\stampatello Theorem}}$ , {\stampatello produce the following } $\underline{\hbox{\stampatello Corollary}}$ : 
$$
C_{f,g_N}(N,1)=C_{f,g}(N,1)
\quad \hbox{\stampatello and} \quad
C_{f,g_N}(N,2)=C_{f,g}(N,2). 
$$

\medskip

\par
\noindent
In all, {\stampatello a T-S Correlation with } $g\neq \0$\enspace {\stampatello and } $g_N\neq \0$ {\stampatello has, from $\UN$ : }
$$
C_{f,g}(N,1)=C_{f,g_N}(N,U_N+1)
\quad \hbox{\stampatello and} \quad
C_{f,g}(N,2)=C_{f,g_N}(N,U_N+2).
\leqno{\CI}
$$

\medskip

\par
\noindent
These are {\stampatello our most interesting Combinatorial Identities}.  
  
\vfill
\eject

\par				
\noindent
In fact, 
$$
C_{f,g}(N,1)=\sum_{2<p\le N}f(p)g^{\ODD}(p+1)
 =\sum_{2<p\le N}f(p)g\left({{p+1}\over {2^{v_2(p+1)}}}\right), 
$$
\par
\noindent
{\stampatello with\enspace $g$\enspace (IPP); and}
$$
C_{f,g}(N,2)=\sum_{2<p\le N}f(p)g^{\ODD}(p+2)
 =\sum_{2<p\le N}f(p)g(p+2), 
$$
\par
\noindent
{\stampatello with\enspace $g$\enspace (IPP), again.} 
\par
On the other side (of {\stampatello \lq \lq huge\rq \rq \enspace shifts), recalling $U_N$ is odd}, 
$$
C_{f,g_N}(N,U_N+1)=\sum_{2<p\le N}f(p)g^{\ODD}_{N}(p+U_N+1)
 =\sum_{2<p\le N}f(p)g_{N}(p+U_N+1), 
$$
\par
\noindent
{\stampatello with\enspace $g_N$\enspace (IPP); while}
$$
C_{f,g_N}(N,U_N+2)=\sum_{2<p\le N}f(p)g^{\ODD}_{N}(p+U_N+2)
 =\sum_{2<p\le N}f(p)g_{N}\left({{p+U_N+2}\over {2^{v_2(p+U_N+2)}}}\right), 
$$
\par
\noindent
again {\stampatello with\enspace $g_N$\enspace (IPP).} 

\bigskip

\par
In passing, previously quoted PNT justifies our terminology about \lq \lq {\stampatello huge shifts}\rq \rq, involving \enspace $U_N$ : 
$$
U_N\EqByDef \prod_{2<p\le N}p = {1\over 2}\exp\left(\sum_{p\le N}\log p\right)
 \EqByDef e^{\theta(N)}/2, 
$$
{\stampatello making } $e^{N(1+o(1))}\ll U_N\ll e^N$, {\stampatello as } $N\to \infty$, {\stampatello from} (quoted above) {\stampatello PNT}.  
  
\vfill
\eject

\par				
\noindent{\bf 5. Main results for Two-Seasons Correlations} 
\bigskip
\par
\noindent
We give a precise and compact exposition of our two results, that we already saw in previous sections, in next $\S5.1$. While we give Main result for the Artifact, starting from definition, in next $\S5.2$. The way we arrive to the Artifact is told in next section, soon after the Proofs of our Main three results in this section.  

\bigskip
\bigskip

\par
\noindent{\bf 5.1. Diophantine Entanglement. Combinatorial Identities} 
\bigskip
\par
\noindent
We start with our newest result: the one we have summarized as \lq \lq {\stampatello Diophantine Entanglement}\rq \rq, above.
\smallskip
\par
\noindent {\bf Theorem 5.1}. {\it Let } $C_{f,g}(N,a)$ {\it be a} {\stampatello T-S Correlation}. {\it Then}
$$
C_{f,g}(N,a)=\sum_{2<p\le N}f(p)g(p+a),
\quad 
\forall a\in \N \enspace \hbox{\stampatello even}
\leqno{\TEQone}
$$
\par
\noindent
{\stampatello and}
$$
C_{f,g}(N,a)=\sum_{2<p\le N}f(p)g((p+a)_{\rm ODD})
 =\sum_{1\le j\le {{\log(N+a)}\over {\log 2}}}\sum_{2<p\le N}f(p)g\left({{p+a}\over {2^j}}\right) ,
\enspace 
\forall a\in \N \enspace \hbox{\stampatello odd};
\leqno{\TEQtwo}
$$
\par
\noindent
{\it which, in case of } {\stampatello Artifact} {\it defined in Theorem $5.3$, i.e. } {\stampatello for next } $f=\mu^2\cdot \1_{\rm ODD}\cdot \Lambda$ {\stampatello and for next } $g=\Lambda_N^{\ODD}$ {\it become}, $\forall N\in \N$ {\it under } {\stampatello Axiom } $4$ ({\it since } {\stampatello Artifact} {\it is a } {\stampatello T-S} : {\it see Th.m } $5.3$) 
$$
C_{f,g}(N,a)=\sum_{2<p\le N}(\log p)\Lambda(p+a),
\quad 
\forall a\in \N \enspace \hbox{\stampatello even}
\leqno{\TEQAone}
$$
$$
C_{f,g}(N,a)=\sum_{1\le j\le {{\log(N+a)}\over {\log 2}}}\sum_{2<p\le N}(\log p)\Lambda\left({{p+a}\over {2^j}}\right),
\quad 
\forall a\in \N \enspace \hbox{\stampatello odd}
\leqno{\TEQAtwo}
$$

\bigskip

\par
\noindent
We expose here the \lq \lq {\stampatello Combinatorial Identities}\rq \rq, coming from periodicity properties above.
\smallskip
\par
\noindent {\bf Theorem 5.2}. {\it Let } $C_{f,g}(N,a)$ {\it be a} {\stampatello T-S Correlation}, {\it with } $g\neq \0$. {\it Recall, from $\S4.1$, the definition of } {\stampatello Universal Period} $U_N$. {\it Then we have the} {\stampatello Combinatorial Identities}
$$
C_{f,g}(N,1)=C_{f,g}(N,U_N+1)
\quad \hbox{\stampatello and} \quad
C_{f,g}(N,2)=C_{f,g}(N,U_N+2). 
\leqno{\TEQthree}
$$

\bigskip
\bigskip

\par
\noindent{\bf 5.2. Main result for Hardy-Littlewood Models: the Artifact} 
\bigskip
\par
\noindent
We give a kind of \lq \lq approximation\rq \rq, for the H-L Correlation $C_{\Lambda,\Lambda}(N,a)$, that is also a T-S Correlation: we call it the {\stampatello Artifact} in the following. 
\smallskip
\par
\noindent {\bf Theorem 5.3}. {\it Define the } {\stampatello Artifact} {\it as the Correlation}
$$
C_{\mu^2\,I_{\rm ODD}\,\Lambda, \Lambda_N^{\rm ODD}}(N,a)\defineq \sum_{2<p\le N}(\log p)\Lambda_N^{\ODD}(p+a),
\quad \forall a\in \N. 
$$ 
\par
\noindent
{\it Then}, $\forall N\in \N$ {\it satisfying} {\stampatello Axiom 4}, {\it it is a} {\stampatello T-S Correlation} : {\it in particular},
$$
a\in \N \enspace \hbox{\stampatello even} \enspace
\enspace \Rightarrow \enspace \enspace 
C_{\mu^2\,I_{\rm ODD}\,\Lambda, \Lambda_N^{\rm ODD}}(N,a)=C_{\mu^2\,I_{\rm ODD}\,\Lambda, \Lambda_N}(N,a),
$$
\par
\noindent
{\it whence, recalling } $L:=\log N$, 
$$
C_{\Lambda,\Lambda}(N,a)=C_{\mu^2\,I_{\rm ODD}\,\Lambda, \Lambda_N^{\rm ODD}}(N,a)+O\left(\left(\sqrt{N}+a\right)L\log(N+a)\right),
\quad \hbox{\stampatello uniformly}\quad \forall a\in \N \enspace \hbox{\stampatello even}. 
\leqno{\TEQfour}
$$

\vfill
\eject

\par				
\noindent{\bf 6. Proofs of Main results above. Hardy-Littlewood Models: towards the Artifact}
\bigskip
\par
\noindent
Since the details of next two Proofs may be found, respectively, in $\S3$ and $\S4$, actually following Proofs are more Sketches of Proofs; but they give a quick glimpse to the ingredients we use. 
\smallskip
\par
The QED stands for the end of a Proof's part; while the Proof ends, as usual, with the symbol\hfill $\square$ 

\bigskip

\par
\noindent {\bf Proof of Theorem 5.1.} We apply the {\stampatello Axioms} \lq \lq everywhere\rq \rq, to prove $\TEQone$ and $\TEQtwo$. 
\par
\noindent
While $\TEQone$ follows from $p+a$ ODD \& $(1.1)$, with trivial $(p+a)_{\rm ODD}=p+a$, our $\TEQtwo$ has $p+a$ EVEN, whence $(1.1)$ uses $(p+a)_{\rm ODD}=(p+a)/v_2(p+a)$, instead.\hfill $\square$ 

\bigskip

\par
\noindent {\bf Proof of Theorem 5.2.} From {\stampatello Axioms}, $g_N=g$, so $g_N\neq \0$ by hypothesis and from $\UN$ we get
$$
\forall a\in \N, \quad C_{f,g}(N,a)=C_{f,g}(N,a+U_N). 
$$
\par
\noindent
From {\stampatello Axioms} and again hypothesis $g\neq \0$, we get $\CI$, in which $g_N=g$ ({\stampatello Axiom $(0)$}, here) entails $\TEQthree$.\hfill $\square$ 
\bigskip

\par
Next Proof is an example of elementary tools in A.N.T., Analytic Number Theory. 
\smallskip
\par
\noindent {\bf Proof of Theorem 5.3.} Recalling von Mangoldt's $\Lambda(n)$ definition soon after $\HLzero$, 
$$
C_{\Lambda,\Lambda}(N,a)\EqByDef \sum_{{n\le N}\atop {n \, \ODD}}\Lambda(n)\Lambda(n+a)+(\log 2)\sum_{K\le {{\log N}\over {\log 2}}}\Lambda(2^K+a)
 =C_{I_{\rm ODD}\,\Lambda,\Lambda}(N,a)+O(L\log(N+a))= 
$$
$$
=C_{I_{\rm ODD}\,\Lambda,\Lambda_N}(N,a)+O(L\log(N+a))-\sum_{N<d\le N+a}\mu(d)(\log d)\sum_{{n\le N}\atop {{n \, \ODD}\atop {n\equiv -a(\bmod d)}}}\Lambda(n)=
$$
$$
=C_{I_{\rm ODD}\,\Lambda,\Lambda_N}(N,a)+O(aL\log(N+a)), 
$$
\par
\noindent
uniformly $\forall a\in \N$, using $\Lambda'(d)=-\mu(d)\log d$, $\forall d\in \N$, see [T]; and since
$$
C_{I_{\rm ODD}\,\Lambda,\Lambda_N}(N,a)=C_{\mu^2\,I_{\rm ODD}\,\Lambda,\Lambda_N}(N,a)+O\left(\left|\sum_{1<K\le {{\log N}\over {\log 2}}}\sum_{2<p\le N^{1/K}}(\log p)\Lambda(p^K+a)\right|\right)=
$$
$$
=C_{\mu^2\,I_{\rm ODD}\,\Lambda,\Lambda_N}(N,a)+O\left(\sum_{1<K\le {{\log N}\over {\log 2}}}N^{1\over K}\log(N+a)\right),
\enspace \forall a\in \N
$$
\par
\noindent
by \v{C}ebi\v{c}ev bound [T], we get 
$$
C_{\Lambda,\Lambda}(N,a)=C_{\mu^2\,I_{\rm ODD}\,\Lambda,\Lambda_N}(N,a)+O\left(\left(\sqrt{N}+a\right)L\log(N+a)\right),
\enspace \forall a\in \N,
$$
\par
\noindent
whence restricting to $a\in \N$ even, we have the approximation $\TEQfour$.\hfill QED
\par
\noindent
Just to complete, we leave the check, of $\BH$ and {\stampatello Axioms} 1,2,3, as an exercise, to the diligent reader.\hfill $\square$ 

\vfill
\eject

\par				
\noindent{\bf 6.1. Hardy-Littlewood Models: towards the Artifact}
\bigskip
\par
\noindent
We will not give explicit remainders, now, when we'll talk about approximating one Correlation with the other, for two reasons: they may be easily obtained by the calculations in previous Proof of Theorem $5.3$ AND they are standard tools in A.N.T. and the use of Theorem 2.1, to cut divisors, from [C1] (often quoted above). See: all of next Correlations $\HLone$ to $\HLfour$ may be considered as \hfill \lq \lq {\stampatello Hardy-Littlewood Models}\rq \rq. 

\bigskip

\par
We start with the {\stampatello Hardy-Littlewood Correlation}, say H-L Corr., i.e. $\HLzero$ above: 
$$
C_{\Lambda,\Lambda}(N,a)=\sum_{n\le N}\Lambda(n)\Lambda(n+a)
\quad
\forall a\in \N.
$$
\par
\noindent
{\stampatello Then, we truncate}, with remainder $O(aL\log(N+a))$, when approximating $\HLzero$:
$$
C_{\Lambda,\Lambda_N}(N,a)=\sum_{n\le N}\Lambda(n)\Lambda_N(n+a)
\quad
\forall a\in \N.
\leqno{\HLone}
$$
\par
\noindent
Now, this is a $\BH-${\stampatello Correlation}; in particular, it {\stampatello has explicit Carmichael-Wintner coefficients}, see Theorem 3.1 [C1]. We don't need to \lq \lq {\stampatello IPP-ify}\rq \rq \thinspace the shift-carrying factor, as \enspace $\Lambda_N$ $\IPP$, already. 
\par
However, we need to \lq \lq {\stampatello ODD-ify}\rq \rq \thinspace our\enspace $\Lambda_N$\enspace into\enspace $\Lambda_N^{\ODD}$\enspace : 
$$
C_{\Lambda,\Lambda_N^{\rm ODD}}(N,a)=\sum_{n\le N}\Lambda(n)\Lambda_N^{\ODD}(n+a)
\quad
\forall a\in \N.
\leqno{\HLtwo}
$$
\par
\noindent
This is \lq \lq {\stampatello close to}\rq \rq (compare Theorem $5.3$ Proof, above) the other
$$
C_{I_{\rm ODD}\Lambda,\Lambda_N^{\rm ODD}}(N,a)=\sum_{{n\le N}\atop {n\,{\rm ODD}}}\Lambda(n)\Lambda_N^{\ODD}(n+a)
\quad
\forall a\in \N.
\leqno{\HLthree}
$$
\par
\noindent
Present Correlation has the \lq \lq {\stampatello Diophantine Entanglement}\rq \rq, which we saw more in general in $\S3$, since : 
$$
C_{I_{\rm ODD}\Lambda,\Lambda_N^{\rm ODD}}(N,a)=\sum_{n\le N}\Lambda(n)\Lambda_N(n+a)
\quad
\forall a\in \N \enspace \hbox{\stampatello even}
\leqno{\HLfour}
$$
\par
\noindent
where RHS (Right Hand Side) is {\stampatello close to H-L Correlation} above; however, the LHS (Left Hand Side) has still the {\stampatello not-square-free} $n$ : these, {\stampatello in presence of}\enspace $\Lambda(n)$,\enspace {\stampatello are negligible}, compare Theorem $5.3$. 

\vfill
\eject

\par				
\noindent{\bf 7. Appendix: Truncated Divisor Sums} 
\bigskip
\par
\noindent
We prove the link \enspace $g'\leftrightarrow \widehat{g}$, if \thinspace $g$ \thinspace is a Truncated Divisor Sum, defined next, with Eratosthenes Transform $g'$. 
\smallskip
\par
\noindent$\underline{\bf Lemma.}$ {\it } {\it Let } $g:\N \rightarrow \C$ {\it be a } {\stampatello Truncated Divisor Sum, by definition} $g=\0$ {\stampatello or, if } $g\neq \0$, {\it we have} {\stampatello by definition}
$$
\exists D\defineq \max(\supporto(g'))\in \N. \quad (\hbox{\stampatello otherwise,}\; \hbox{\it in case}\enspace g=\0,\enspace \hbox{\stampatello we may set}\enspace D\defineq 0.)
$$
\par
\noindent
{\it Then, see } [C0], $\exists \WinT g$, $\exists \CarT g$ {\it and the } {\stampatello unique Ramanujan Coefficient} {\it of } $g$ {\it is } $\widehat{g}\defineq \WinT g=\CarT g$, 
$$
\forall q\in\N, \enspace \widehat{g}(q)=\Win_q\; g\defineq \sum_{{d\le D}\atop {d\equiv 0(\!\!\bmod q)}}{{g'(d)}\over d} \enspace 
\hbox{\stampatello is the $q-$th Ramanujan Coefficient of}\enspace g, 
\leqno{(\ast)}
$$
\par
\noindent
{\stampatello with} : $D=\max(\supporto(\widehat{g}))\in \N$, $\forall g\neq \0$. ($g=\0$ {\it has} \thinspace $D=0$: $\widehat{\0}=\WinT \0=\CarT \0=\0$.) {\it These give, $\forall g\neq \0$, a} {\stampatello fixed-length Ramanujan Expansion.} ({\it See Proof.}) {\stampatello Furthermore,} {\it in same hypotheses \& notations,} 
$$
\hbox{\stampatello if } \enspace \exists \hbox{\script S}\subseteq \N, \hbox{\stampatello divisor-closed, with } \enspace
\supporto(g')\subseteq \hbox{\script S}, 
\enspace \hbox{\stampatello then } \enspace 
\supporto(\WinT g)\subseteq \hbox{\script S} 
\leqno{(\ast\ast)}
$$
\par
\noindent
{\stampatello and}
$$
\hbox{\stampatello if } \enspace \exists \hbox{\script S}\subseteq \N, \hbox{\stampatello divisor-closed, with } \enspace
\supporto(\WinT g)\subseteq \hbox{\script S},  
\enspace \hbox{\stampatello then } \enspace 
\supporto(g')\subseteq \hbox{\script S}. 
\leqno{(\ast\ast\ast)}
$$
\smallskip
\par
\noindent$\underline{\bf Proof.}$ By $(1.8)$ [C1], say \lq \lq {\stampatello Ramanujan Orthogonality}\rq \rq, abbrev. R.O.: 
${\displaystyle 
\1_{d|a}={1\over d}\sum_{q|d}c_q(a),
\enspace 
\forall a,d\in \N
}$,  
$$
\forall a\in \N, 
\enspace 
g(a)=\sum_{{d\le D}\atop {d|a}}g'(d)
 \buildrel{R.O.}\over{=\!=\!=} \sum_{d\le D}{{g'(d)}\over d}\sum_{q|d}c_q(a)
  = \sum_{q\le D}\sum_{{d\le D}\atop {d\equiv 0(\!\!\bmod q)}}{{g'(d)}\over d}c_q(a)
   \buildrel{(\ast)}\over{=\!=\!=}\sum_{q\le D}(\Win_q\; g)c_q(a); 
$$
\par
\noindent
namely we get, having assumed $g\neq \0$, like henceforth of course, the 
$$
\forall a\in \N, 
\enspace 
g(a)=\sum_{q\le D}(\Win_q\; g)c_q(a)
     =\sum_{q\le D}\widehat{g}(q)c_q(a).
\leqno{\RamaExp} 
$$
\par
\noindent
We get $D=\max(\supporto(\widehat{g}))$, as: $\widehat{g}(q)=0$, $\forall q>D$ by $(\ast)$ above and by $D$ definition: $\widehat{g}(D)\EqByDef {{g'(D)}\over D}\neq 0.
$
\QED
\par
\noindent
See that, in $\RamaExp$ {\stampatello of $g$, above, $\Win_q\; g$ doesn't depend on $a$, so:} 
$$
\forall d\in \N, 
\enspace 
g'(d)=d\sum_{{K\le D/d}\atop {dK\in \supporto(\WinT g)}}\mu(K)\Win_{dK}\; g 
\leqno{\LuchtExp} 
$$
\par
\noindent
holds, too, writing \enspace $c_q(a)={\displaystyle \sum_{{d|q}\atop {d|a}} }d\mu\left(q\over d\right)$ \enspace (see [C1], $\S1$) and setting $K:=q/d$ from {\stampatello M\"obius Inversion in Ramanujan Expansion} above. \hfill {\stampatello From $(\ast)$ above, we get $(\ast\ast)$,} as: if \thinspace $\hbox{\script S}\subseteq \N$ \thinspace is {\stampatello divisor-closed}, then 
\smallskip
\centerline{
$
\forall q\in\N , \enspace \Win_q\; g\EqByDef {\displaystyle \sum_{{d\le D}\atop {{d\equiv 0(\!\!\bmod q)}\atop {d\in \supporto(g')}}} } {{g'(d)}\over d}
\enspace \hbox{\stampatello entails}, \hbox{\rm when } \supporto(g')\subseteq \hbox{\script S}, \hbox{\rm that}: 
\enspace q\in \supporto(\WinT g) \thinspace \Rightarrow \thinspace q\in \hbox{\script S}.
$
\hfill QED}
\par
\noindent
{\stampatello From $\LuchtExp$ above, we prove } $(\ast\ast\ast)$, since: if \thinspace $\hbox{\script S}\subseteq \N$ \thinspace is {\stampatello divisor-closed}, then 
$$
\supporto(\WinT g)\subseteq \hbox{\script S}\enspace \hbox{\stampatello entails}: 
\enspace d\in \supporto(g') \thinspace \Rightarrow \thinspace d\in \hbox{\script S}.
$$
\hfill $\square$ 
 
\vfill
\eject

\par				
\noindent
We recall that any \enspace $g:\N \rightarrow \C$, $\forall N\in \N$ fixed, has {\stampatello $N-$truncation} : 
$$
g_N(m)\defineq \sum_{{d|m}\atop {d\le N}}g'(d),
\quad
\hbox{\stampatello where}\enspace \hbox{\rm [T]} 
\enspace g'\defineq \mu \ast g \enspace
\hbox{\stampatello is the Eratosthenes Transform of our $g$.}
$$
\par
\noindent
Here {\stampatello we may have, in general, $g'(N)=0$, too;} i.e., $N$ above is a {\stampatello Wintner's Range, but we are not assuming $N$ is the exact Wintner's Range} (see page 7 in [C2]: this happens exactly IFF \enspace $g'(N)\neq 0$). 
\par
\noindent
Recall the definition of {\stampatello Correlation} of $f,g:\N \rightarrow \C$ (see [C1], $\S2.2$) as \enspace $C_{f,g}(N,a)\defineq {\displaystyle \sum_{n\le N}} f(n)g(n+a)$, where $a\in \N$ is the {\stampatello shift} (the variable of our Correlation, as a kind of arithmetic function), while $N\in \N$ is the {\stampatello length}, a kind of \lq \lq large\rq \rq, fixed integer.
\par
From this definition \& the above $N-$truncation, we get (compare [C1], Theorem 2.1) the following.
\smallskip 
\par
\noindent$\underline{\bf Proposition.}$ {\it Let } $f,g:\N \rightarrow \C$ {\it be } {\stampatello two arbitrary arithmetic functions.} {\it Then,} \enspace $\forall a,N\in \N$, 
$$
C_{f,g}(N,a)-C_{f,g_N}(N,a)=\sum_{N<d\le N+a}g'(d)\sum_{{n\le N}\atop {n+a\equiv 0(\!\!\bmod d)}}f(n). 
$$
\medskip
\par
The $\underline{\hbox{\rm \lq \lq small\rq \rq}}$ values of $\underline{\hbox{shifts}}$ $a\in \N$ $\underline{\hbox{\rm w.r.t.}}$ the $\underline{\hbox{\rm length}}$ $N\in \N$ $\underline{\hbox{\rm give}}$ the $\underline{\hbox{\rm following}}$.
\smallskip 
\par
\noindent$\underline{\bf Theorem.}$ {\it Fix any } $f,g:\N \rightarrow \C$. {\it Assume } $a,N\in \N$ {\it with } $a\le N$. {\it Then, from above Proposition,}
$$
C_{f,g}(N,a)-C_{f,g_N}(N,a)=\sum_{N<d\le N+a}g'(d)f(d-a), 
\enspace \hbox{\it whence } \enspace 
C_{f,g}(N,1)-C_{f,g_N}(N,1)=g'(N+1)f(N). 
$$
\par
\noindent
{\it Furthermore, in the same hypotheses above, provided $N>1$, we get}    
$$
C_{f,g}(N,2)-C_{f,g_N}(N,2)=g'(N+1)f(N-1)+g'(N+2)f(N). 
$$
\smallskip
\par
\noindent$\underline{\bf Proof.}$ We've: $n+a\equiv 0(\bmod \; d)\Leftrightarrow\exists K\in \N: n+a=Kd$; $N<d,n\le N$ $\Rightarrow$ $K<1+{a\over N}\le 2$ $\Rightarrow$ $K=1$.\hfill $\square$ 
\medskip
\par				
\centerline{\stampatello Bibliography}
\medskip
\item{[C0]} G. Coppola, {\sl Recent results on Ramanujan expansions with applications to correlations}, Rend. Sem. Mat. Univ. Pol. Torino {\bf 78.1} (2020), 57--82. 
\smallskip
\item{[C1]} G. Coppola, {\sl General elementary methods meeting elementary properties \thinspace of \thinspace correlations}, {\tt available}\break{\tt online at} \enspace arXiv:2309.17101 (2nd version)
\smallskip
\item{[C2]} G. Coppola, {\sl On Ramanujan smooth expansions for a general arithmetic function}, {\tt available online at} \enspace arXiv:2407.19759v1 
\smallskip
\item{[C3]} G. Coppola, {\sl Good Ramanujan Expansions: A suitably enhanced decay of coefficients has important consequences}, {\tt available online at} \enspace arXiv:2509.24456v1 
\smallskip
\item{[CM]} G. Coppola and M. Ram Murty, {\sl Finite Ramanujan expansions and shifted convolution sums of arithmetical functions, II}, J. Number Theory {\bf 185} (2018), 16--47. 
\smallskip
\item{[HL]} G.H. Hardy and J.E. Littlewood, {\sl SOME PROBLEMS OF 'PARTITIO NUMERORUM'; III: ON THE EXPRESSION OF A NUMBER AS A SUM OF PRIMES.} Acta Mathematica {\bf 44} (1923), 1--70. 
\smallskip
\item{[M]} M. R. Murty, {\sl Ramanujan series for arithmetical functions}, Hardy-Ramanujan J. {\bf 36} (2013), 21--33. Available online 
\smallskip
\item{[R]} S. Ramanujan, {\sl On certain trigonometrical sums and their application to the theory of numbers}, Transactions Cambr. Phil. Soc. {\bf 22} (1918), 259--276.
\smallskip
\item{[T]} G. Tenenbaum, {\sl Introduction to Analytic and Probabilistic Number Theory}, Cambridge Studies in Advanced Mathematics, {46}, Cambridge University Press, 1995. 
\medskip
\par
\leftline{\tt Giovanni Coppola - Universit\`{a} degli Studi di Salerno (affiliation)}
\leftline{\tt Home address : Via Partenio 12 - 83100, Avellino (AV) - ITALY}
\leftline{\tt e-mail : giocop70@gmail.com}
\leftline{\tt e-page : www.giovannicoppola.name}
\leftline{\tt e-site : www.researchgate.net}

\bye